\documentclass[oneside]{amsart}
\usepackage[mathscr]{eucal}
\usepackage{amssymb, amsmath, array, hyperref}

\newfont{\sheaf}{eusm10 scaled\magstep1}

\newtheorem{theorem}{Theorem}[section]

\newtheorem{THM}{Theorem}
\newtheorem{Q}{Question}
\newtheorem{COR}{Corollary}

\newtheorem{prop}[theorem]{Proposition}
\newtheorem{lemma}[theorem]{Lemma}

\newtheorem{remark}[theorem]{Remark}
\newtheorem{cor}[theorem]{Corollary}

\numberwithin{equation}{section}

\newtheorem{example}[theorem]{Example}

\begin{document}

\title[]
{Stability of Holomorphic Foliations with Split Tangent Sheaf}

\author{Fernando Cukierman}

\author{Jorge Vit\'orio Pereira}

\dedicatory{}

\begin{abstract}
We show that the set of  singular holomorphic  foliations on
projective spaces with split tangent sheaf and  {\it good} singular
set is open in the space of holomorphic foliations. We also give a
cohomological criterion for the rigidity of holomorphic foliations
induced by group actions and prove the existence of rigid
codimension one foliations of degree $n-1$ on $\mathbb P^n$ for
every $n\ge 3$.
\end{abstract}

\thanks{The second  author is supported by Instituto Unibanco and Cnpq. Both  authors were partially supported by  CAPES-SPU}



\maketitle

\tableofcontents

\section{{Introduction  and Statement of Results}}

Our main object of study in this article is the geometry of the
spaces of singular holomorphic distributions  and singular
holomorphic foliations on  complex projective spaces.

Loosely speaking, a codimension $q$ singular holomorphic
distribution on $\mathbb P^n$ is a holomorphic field of
$(n-q)$-planes on the complement of a Zariski closed set  of
codimension at least two. When this  plane field is involutive  we
have a singular holomorphic foliation. The most basic projective
invariant that one can attach to a distribution or to a foliation is
its degree. The degree  is defined as the degree of the tangency of
the distribution or foliation  with a generic $\mathbb P^q$ linearly
embedded in $\mathbb P^n$.

In this work we will denote  by $\mathscr D_q(n,d)$ and $\mathscr
F_q(n,d)$ the spaces of distributions and of foliations on $\mathbb
P^n$ of codimension $q$ and degree $d$. These spaces turn out to be
quasi-projective varieties and  we will give new information about
the irreducible components of $\mathscr F_q(n,d)$ and $\mathscr
D_q(n,d)$ for arbitrary $n\ge 3$ and $1 \le q \le n-1$.

Precise definitions of the relevant  concepts will be given latter
in the Introduction. Before  that we would like to recall some known
results and share our main motivation with the reader.

\subsection{Known Results}

 The systematic study of the
irreducible components of $\mathscr F_1(n,d)$ seems to have been
initiated by Jouanolou in \cite{Jouanolou} where the irreducible
components of $\mathscr F_1(n,0)$ and $\mathscr F_1(n,1)$ where
classified for all $n \ge 3$. The classification of the irreducible
components of $\mathscr F_1(n, 2)$ was achieved by Cerveau and Lins
Neto in \cite{CL}. Besides the classification in low degrees, a few
infinite families of irreducible components of the space of
codimension one foliations on projective spaces are known:
\smallskip

\begin{tabular}{ll}
  (a)   rational
 \cite{GML}; &  (b)  logarithmic
\cite{Omegar};\\  (c)   linear pull-back  \cite{CaLN}; & (d)
  generic pull-back  \cite{CLE}; \\ (e)
associated to  the affine Lie algebra \cite{CCGL}.
\end{tabular}

\medskip

Unlike the codimension one case, where there is a growing
literature, we are aware of just one result about the irreducible
components of $\mathscr F_q(n, d)$ when $q\ge 2$:  the
classification of the irreducible components of $\mathscr F_q(n, 0)$
given in \cite[Proposition 3.1]{CeDe}. We point out that for $d\ge
0$ and $q\ge 2$ no irreducible components of $\mathscr F_q(n,d)$
were known so far. Althought \cite[Theorem A]{Airton} solves a
similar problem for degree one codimension $q$ foliations on
$\mathbb C^n$, i.e., codimension $q$ foliations on $\mathbb C^n$
induced by linear $q$-forms.

\subsection{Motivation: Foliations induced by Group
Actions}\label{S:motivation}

The key question behind the developments here presented was
motivated by a conjecture of Cerveau and Deserti made  in the recent
monograph \cite{CeDe}. There, after classifying the codimension one
foliations of degree $3$ on $\mathbb P^4$ induced by Lie subalgebras
of $\mathfrak{aut}(\mathbb P^4)\cong \mathfrak{sl}(5,\mathbb C)$,
they conjecture that one of these foliations, more precisely the one
that admits
\[
\frac{ \left( z_0 z_4^3 - (2z_1z_3 + z_2^2)z_4^2 + 2z_2z_3^2z_4 -
\left(\frac{ z_3^4}{2}\right) \right)^3 } { \left( z_1z_4^2 -
z_2z_3z_4 + \frac{z_3^3}{3} \right)^4 } \,
\]
as rational first integral is rigid, that is, there exists an
irreducible component of $\mathscr F_1(4,3)$ whose generic element
is projectively equivalent to the one induced by the levels of the
rational function above ({\it cf.}  introduction of {\it loc.cit.}).
With this conjecture in mind we were
naturally lead to the following:

\begin{Q}\label{Q:1}
Under which conditions a Lie subalgebra
 $\mathfrak{g}$ of $\mathfrak{sl}(n+1,\mathbb C)$ induces a rigid foliation of
$\mathbb P^n$?
\end{Q}

If we assume that $\mathcal F$ is a foliation induced by a
subalgebra  $\mathfrak{g} \subset \mathfrak{aut}(\mathbb P^n) \cong
\mathfrak{sl}(n+1,\mathbb C)$ and that the integration of $\mathfrak
g$ induces an action which is locally free on the complement of an
algebraic subvariety $\Sigma \subset \mathbb P^n$ of codimension at
least two then the tangent sheaf of $\mathcal F$ is trivial. Thus,
at least on this case, we see that Question \ref{Q:1} is strictly
related to

\begin{Q}\label{Q:2}
Under which conditions a deformation of a foliation $\mathcal F$
with trivial tangent sheaf still has trivial tangent sheaf?
\end{Q}

After a careful study  of the myriad of examples presented in
\cite{CeDe} we realized that  good candidates for sufficient
conditions for codimension $q$ foliations in Question \ref{Q:2} are:
\begin{equation}\label{E:condition}
\left\lbrace
\begin{array}{clc}
\mathrm{codim} \ \mathrm{sing}(d\omega) \ge 3 & \text{ when } & q=1\,  \\
\mathrm{codim} \ \mathrm{sing}(\omega) \ge 3 & \text{ when } & q\ge
2 \,
\end{array}
\right. \,
\end{equation}
where $\omega$ is a homogeneous $q$-form defining $\mathcal F$.

\medskip

At this point, in order to keep the prose intelligible,    it is
clear that we need to be more precise about some basic notions.

\subsection{Basic Definitions}
A {\bf singular holomorphic distribution} $\mathcal D$ of dimension
$p$ on a projective space   $\mathbb P^n$ is a rational section  of
$\mathbb G_p(\mathrm T\mathbb P^n)$, where $\mathbb G_p(\mathrm
T\mathbb P^n)$ is the Grassmann bundle of $p$-planes in $\mathrm
T\mathbb P^n$. The distribution  $\mathcal D$ can also be dually
presented as a rational section of $\mathbb G_q(\mathrm T^*\mathbb
P^n)$, where $q$ is the codimension of $\mathcal D$, {\it i.e.},
$p+q=n$.

If we consider the standard embedding of $\mathbb G_q(\mathrm
T^*\mathbb P^n)$ in $\mathbb P(\Omega^q_{\mathbb P^n})$ then the
rational section defining $\mathcal D$ can be interpreted as the
projectivization  of the image of a rational $q$-form $ \omega$. If
$(\omega)_0$  is the divisorial part of the zero scheme of $\omega$,
$(\omega)_{\infty}$ is the divisor of poles of $\omega$ and  we set
\[
\mathcal L = \mathcal O_{\mathbb P^n}\big( ( \omega)_{\infty} - (
\omega)_{0} \big)
\]
then the rational section defining $\mathcal D$ is the
projectivization of the image of some $\omega \in \mathrm
H^0(\mathbb P^n,\Omega^q_{\mathbb P^n} \otimes \mathcal L)$ with
zero (or singular)  set of codimension at least two.

The {\bf singular set} of $\mathcal D$, denoted by
$\mathrm{sing}(\mathcal D)$, is the zero set of the twisted $q$-form
$\omega$. Notice that by definition $\mathrm{sing}(\mathcal D)$ has
always codimension at least two. The {\bf degree} of $\mathcal D$,
denoted by $\deg(\mathcal D)$, is by definition the degree of the
zero locus of $i^*\omega$, where $i:\mathbb P^q \to \mathbb P^n$ is
a linear embedding of a generic $q$-plane. Since $\Omega^q_{\mathbb
P^q} = \mathcal O_{\mathbb P^q}(-q-1)$ it follows at once that
$\mathcal L = \mathcal O_{\mathbb P^n}(\deg(\mathcal D)+ q + 1)$.

\medskip
If  $\omega \in \mathrm H^0(\mathbb P^n, \Omega^q_{\mathbb
P^n}\otimes \mathcal L)$ and $q>1$ then in general the
projectivization of the graph of $\omega$ is not contained in
$\mathbb G_q(\mathrm T^*\mathbb P^n)$. It will be the case ({\it
cf.} for instance \cite{Griffiths-Harris,Airton}) if, and only if,
$\omega$ satisfies (pointwise) the Pl\"{u}cker conditions.

It is well known that the vector space $\mathrm H^0(\mathbb P^n,
\Omega^q_{\mathbb P^n}\otimes \mathcal L)$ can be canonically
 identified with the vector space of $q$-forms on $\mathbb
C^{n+1}$ with homogeneous coefficients of degree $d+1$ whose
contraction with
 the radial (or Euler) vector field $R=\sum_{i=0}^n x_i
\frac{\partial}{\partial x_i}$ is identically zero, {\it cf.}
\cite{Jouanolou}. Taking advantage of this identification   the
Pl\"{u}cker equations can be written on $\mathbb C^{n+1}$ as
\begin{equation}\label{E:pluckeri}
 ( {i}_v\omega) \wedge \omega    = 0 \quad \text{ for every }v \in \bigwedge^{q-1}\mathbb
 C^{n+1}.
\end{equation}

When  $\omega \in \mathrm H^0(\mathbb P^n, \Omega^q_{\mathbb
P^n}(d+q+1)$ satisfies (\ref{E:pluckeri}) then the kernel of the
morphism of sheaves
 \begin{eqnarray*}
 \mathrm T{\mathbb P^n} &\longrightarrow&
\Omega^{q-1}_{\mathbb P^n}(d+q+1) \\
v &\mapsto& i_v \omega
\end{eqnarray*}
 defined by contraction with $\omega$ has
generic rank $q$ and is a  sheaf  called the {\bf tangent
sheaf} of $\mathcal D$, denoted in this work by $\mathrm T \mathcal
D$. Alternatively one could define a codimension $q$ distribution on
$\mathbb P^n$ as a rank $n-q$ saturated subsheaf of $\mathrm T\mathbb P^n$,
where saturated means that the corresponding cokernel is  torsion-free. Both
definitions turn out to be equivalent and for our purposes the
definition in terms of twisted differential forms satisfying Pl\"{u}cker
equations will be more manageable.

\medskip

A {\bf codimension $q$ singular holomorphic foliation} $\mathcal F$
on  $\mathbb P^n$ is a codimension $q$ distribution $\mathcal F$
with tangent sheaf closed under Lie bracket, i.e., $[\mathrm T
\mathcal F, \mathrm T\mathcal F] \subset \mathrm  T \mathcal F$. If
$\omega \in \mathrm H^0(\mathbb P^n, \Omega^q_{\mathbb P^n}\otimes
\mathcal L)$ is a $q$-form defining $\mathcal F$ then $\omega$
satisfies (\ref{E:pluckeri}) and the involutiveness of  $\mathrm
T\mathcal F$ is equivalent to, {\it cf.} \cite[proposition
1.2.2]{Airton},
\begin{equation}\label{E:pluckerii}
 ( {i}_v\omega) \wedge d\omega    = 0 \quad \text{ for every }v \in \bigwedge^{q-1}\mathbb C^{n+1}   .
\end{equation}
Of course $\omega$ and $d\omega$ on the formula above are
homogeneous differential forms on $\mathbb C^{n+1}$.

\medskip

We will denote $\mathscr D_q(n,d) $ (resp. $\mathscr F_q(n,d)$)  the
quasi-projective subvariety of $\mathbb P\mathrm H^0(\mathbb P^n,
\Omega^q_{\mathbb P^n}(d+q+1))$ whose points parametrize the degree
$d$ and codimension $q$ distributions (resp. foliations) on $\mathbb
P^n$, i.e.,
\[
\begin{array}{lclcl}
  \mathscr D_q(n,d) &=& \{ [\omega]    \quad |  \quad
  \omega  \text{ satisfies } (\ref{E:pluckeri})   & \text{ and }&   \mathrm{codim} \
\mathrm{sing}(\omega) \ge 2 \} \, ; \\
  \mathscr F_q(n,d) &=& \lbrace [\omega]    \quad  |
  \quad
  \omega  \text{ satisfies } (\ref{E:pluckeri}),(\ref{E:pluckerii}) & \text{ and }&   \mathrm{codim} \
\mathrm{sing}(\omega) \ge 2 \rbrace .
\end{array}
\]

In words: $\mathscr D_q(n,d)$ (resp. $\mathscr F_q(n,d)$) is the
space of codimension $q$ singular holomorphic  distributions (resp.
foliations) of degree $d$ on $\mathbb P^n$.

\subsection{Main Results}
Our first main result says that the conditions (\ref{E:condition})
 are indeed sufficient  for the stability of trivial
tangent sheaf for codimension one foliations. In fact we are able to
settle the more general:

\begin{THM}\label{T:splits1}
Let $n\ge 3$, $d\ge 0$ be integers and $\mathcal F = [\omega] \in
\mathscr F_1(n,d)$ be a singular holomorphic  foliation on $\mathbb
P^n$. If $\mathrm{codim} \ \mathrm{sing}(d\omega) \ge 3$ and
\[
T\mathcal F \cong \bigoplus_{i=1}^{n-1} \mathcal O_{\mathbb
P^n}(e_i), \quad e_i \in \mathbb Z,
\]  then there exists a Zariski-open neighborhood
$\mathcal U \subset \mathscr F_1(n,d)$ of $\mathcal F$ such that
$T\mathcal F'\cong \oplus_{i=1}^{n-1} \mathcal O_{\mathbb P^n}(e_i)$
for every $\mathcal F' \in \mathcal U$.
\end{THM}

In dimension $3$ a variant of  the above Theorem appears as the
first step of the proof of \cite[Theorem 1]{CCGL}. There the
argumentation is based on  the deformation
theory of holomorphic vector bundles and a detailed understanding of
germs of integrable $1$-forms $\omega$ on $\mathbb C^3$ such that
$d\omega$ has an isolated singularity. Unfortunately such a detailed
understanding is not available in higher dimensions or codimensions.
Our proof, based on infinitesimal techniques, is completely
different and works uniformly in all dimensions.

Our method also works  for codimension $q$ distributions when $q\ge
2$,  and yields the following result.

\begin{THM}\label{T:splits2}
Let $n\ge 4$, $q\ge 2$, $d\ge0$ be integers and  $\mathcal D\in
\mathcal D_q(n,d)$ be a singular holomorphic distribution on
$\mathbb P^n$. If $ \mathrm{codim} \ \mathrm{sing}(\mathcal D) \ge
3$ and
\[
T\mathcal D  \cong \bigoplus_{i=1}^{n-q} \mathcal O_{\mathbb
P^n}(e_i), \quad e_i \in \mathbb Z \, ,
\]
then  there exists a Zariski-open neighborhood $\mathcal U \subset
\mathscr D_q(n,d)$ of $\mathcal D$ such that $T\mathcal D'\cong
\oplus_{i=1}^{n-q} \mathcal O_{\mathbb P^n}(e_i)$ for every
$\mathcal D' \in \mathcal U$.
\end{THM}

\medskip

In \S\ref{S:applications} we present two  immediate consequences
 of Theorems \ref{T:splits1} and \ref{T:splits2}. The
first one is  a generalization of a well-known result of Camacho and
Lins Neto about the linear pull-back of foliations. The second one
is a generalization of a result of Calvo-Andrade, Cerveau, Lins Neto
and Giraldo about the stability of foliations associated to affine
Lie algebras.

\medskip

Combining Theorems \ref{T:splits1} and \ref{T:splits2} with
Richardson's results about  the rigidity of subalgebras of Lie
algebras we are able to give an answer to Question \ref{Q:1}.  It
comes in the form of our third main result.

\begin{THM}\label{T:3}
Let $\mathcal F$ be a codimension $q$ foliation on $\mathbb P^n$
induced by a Lie subalgebra $\mathfrak g \subset \mathfrak{sl}(n+1,
\mathbb C)$ whose corresponding action is locally free outside a
codimension $2$ analytic subset of $\mathbb P^n$. Suppose that
$\mathcal F$ satisfies the hypothesis of Theorem \ref{T:splits1}
(for $q=1$) or Theorem \ref{T:splits2} (for $q \ge 2)$. If $\
\mathrm{H}^1\left(\mathfrak{g}, \ \mathfrak {sl}(n+1, \mathbb C)/
\mathfrak g\right) = 0 \ $ then $\mathcal F$ is rigid, i.e., the
closure of the orbit of such foliation under the automorphism group
of $\mathbb P^n$ is an irreducible component of $\mathscr
F_q(n,n-q)$.
\end{THM}

\medskip

On the one hand, Theorem \ref{T:3} together with  well-known
vanishing results for the cohomology of semi-simple Lie algebras
yields  a handful of new rigid foliations of codimension at least
two on projective spaces. On the other hand  these general vanishing
results are ineffective in the codimension one case: we prove in
proposition \ref{P:contra} that all codimension one foliations
induced by generically locally free actions of semi-simple groups
are not rigid.  Albeit  we are able to prove
that for every $n\ge3$ there exists a codimension one {rigid}
foliation on $\mathbb P^n$ of degree $n-1$ induced by a meta-abelian
subalgebra of $\mathfrak{aut}(\mathbb P^n)$.

\begin{THM}\label{T:infinito}
For every $n\ge 3$ the codimension one foliation on $\mathbb P^n$
induced by the Lie subalgebra of $\mathfrak{aut}(\mathbb P^n)$
generated by
\begin{eqnarray*}
X = \sum_{i=0}^n (n-2i) z_i \frac{\partial z}{\partial z_i} \, \quad
\text{and} \quad  Y_k = \sum_{i=0}^{n-k} z_{i+k} \frac{\partial
z}{\partial z_i}\, ,\quad k=1\ldots n-2 \,
\end{eqnarray*}
is rigid.
\end{THM}

We point out that when $n=3$ the  rigidity of the foliation in
Theorem \ref{T:infinito}  was established in \cite{CL}. For $n=4$
the Theorem gives a positive answer to the Cerveau-Deserti
conjecture mentioned in \S\ref{S:motivation}.

We also found two other examples of  rigid codimension one
foliations induced by group actions: one in $\mathbb P^6$ and the
other in $\mathbb P^7$, {\it cf.} Table 1 below. In contrast with
the Lie algebras presented in Theorem \ref{T:infinito} for these
examples the first derived algebras are not abelian but just
nilpotent.

\medskip

\subsection{New Irreducible Components}

Throughout   the text the reader will find several new irreducible
components of the spaces of foliations. We summarize in the table
below all the rigid foliations  associated to Lie subalgebras of
$\mathfrak{aut}(\mathbb P^n)$ that appear in this work.

All the omitted Lie brackets that can't be  deduced by anti-symmetry from 
the specified ones are understood to be zero. We emphasize that Corollary \ref{C:PB1}
implies that the foliations obtained from the ones in the Table by a
generic linear pull-back  are also rigid.

\begin{table}[ht]
  \centering
\begin{center}
\begin{tabular}{|c|c|c|p{2.3in}|c|}
\hline $\mathbf{q}$  & $\mathbf{n}$ & \bf{degree}  &   \bf{Lie Algebra} & \bf{Reference}   \\
\hline \hline
  $q\ge 1$  & $2+q$ & $2$ &  $\mathfrak{aff}(\mathbb C)$  &  Ex. \ref{E:AFF} \\
\hline
  $ q \ge 2$ & $3+q$ & $3$  &  $\mathfrak{sl}(2,\mathbb C)$ & Ex. \ref{E:PSL} \\
\hline
  $ q \ge 2$ & $n \gg 0$ & $n - q $  &  \text{semi-simple} & Ex.  \ref{E:semisimples} \\
\hline
  $1$ & $n\ge 3$ & $n-1$  & \small{$<X,Y_1,\ldots, Y_{n-2}>$ satisfying $[X,Y_i]=-2iY_i$}   &  Thm  \ref{T:infinito}  \\
\hline
  $1$ & $6$ & $5$  & \small{$<X,Y_1,\ldots, Y_4>$ satisfying $[X,Y_i]=-2iY_i$, $[Y_1,Y_j]=Y_{j+1}$ }   &  Prop.   \ref{P:67}  \\
\hline
  $1$ & $7$ & $6$  & \small{$<X,Y_1,\ldots, Y_5>$ satisfying $[X,Y_i]=-2iY_i$, $[Y_1,Y_j]=Y_{j+1}$, $[Y_2,Y_3] = -\frac{5}{2} Y_5$ }   &  Prop.  \ref{P:67} \\
\hline
\end{tabular}
\end{center}
  \caption{Rigid Foliations in $\mathscr F_q(n,d)$  }\label{tabela}
\end{table}

Besides the rigid foliations and associated irreducible components
of $\mathscr F_q(n,d)$ presented above we  found an infinite family
of irreducible components whose generic element is the linear
pull-back of foliations by curves ({\it cf.} Example
\ref{E:pullback}) and an infinite family of irreducible components
whose generic element is a foliation induced by an abelian action
({\it cf.} Example \ref{E:diagonal}). Both families generalize to
arbitrary codimension previously known examples of codimension one.

\medskip

\section{{Preliminaries} }\label{S:inicial}

\subsection{Calculus on $\mathbb C^{n+1}$}\label{S:calculus}

For vector fields $X, Y$ on $\mathbb C^{n+1}$ we recall that the
interior product and Lie derivative satisfy the following useful
relations
\begin{equation}\label{E:util}
 [ L_X, i_Y ] = i_{[X,Y]};    \, \,  \, \, \, \, \, [ L_X, L_Y ]
= L_{[X,Y]}; \, \, \, \, \,  \, \,  L_X \Omega = \mathrm{div}(X)
\Omega;
\end{equation}
where $\Omega$ denotes the euclidean volume form in $\mathbb
C^{n+1}$, i.e., $\Omega=dx_0 \wedge \cdots \wedge dx_n$. For
instance they imply that
\begin{equation}\label{E:div0}
\mathrm{div}([X,Y]) = X(\mathrm{div}(Y)) - Y(\mathrm{div}(X)) \, .
\end{equation}

If  $\omega$ denotes  a degree $d$ homogeneous $p$-form, i.e. the
coefficients of $\omega$ are homogeneous polynomials of degree $d$,
then
\[
   L_R \omega = (d+p)  \omega \, .
\]
In particular if $\omega$ is annihilated by the radial vector
field  then
\begin{equation}\label{E:dif2}
i_R d \omega = (d+p) \omega \, .
\end{equation}

If $X $ is a degree $d$ homogeneous vector field then
\[
    [X,R] = (1 - d) X \, .
\]

\medskip

The next lemma is a kind of dual version of formula (\ref{E:dif2})
for integrable distributions and will be used in the proof of
Theorem \ref{T:splits1}.

\begin{lemma}\label{L:prepara}
Let $X_1, \ldots, X_q$ be homogeneous polynomial vector fields on
$\mathbb C^{n+1}$ such that the $(n-q)$-form $ \eta = i_{X_1} \cdots
i_{X_q} i_R \Omega $ is integrable and has singular set of
codimension $\ge 2$. Then there exist homogeneous polynomial vector
fields $\widetilde X_1, \ldots, \widetilde X_q$ such that
\begin{enumerate}
\item $\eta = i_{\widetilde X_1} \cdots i_{\widetilde X_q} i_R
\Omega \, ;$
\item $\displaystyle{d\eta = (-1)^q} \left(n +1 - q +\sum_i \deg(X_i)  \right) i_{\widetilde X_1} \cdots i_{\widetilde X_q } \Omega
\, ;$ \item $\deg(X_i) = \deg(\widetilde X_i ) \text{ for every } i
\, .$
\end{enumerate}
\end{lemma}
\begin{proof}
 Let  $ X_{q+1}=R$. Since $\eta$ is integrable,
\[
[X_i,X_j]=\sum_{l=1}^{q+1} a_{ij}^l X_l \,
\]
for some rational functions $a_{ij}^l$. Under our hypothesis  the
rational functions $a_{ij}^l$ are regular everywhere, i.e.,
homogeneous polynomials. In fact
\[
   [X_i, X_j ] \wedge X_1 \cdots \wedge \widehat{X_l} \wedge
   \cdots \wedge X_q \wedge R = \pm a_{ij}^l  X_1 \cdots \wedge {X_l} \wedge
   \cdots \wedge X_q \wedge R \, ,
\]
and since, by hypothesis, the zero set of $X_1 \cdots \wedge {X_l}
\wedge  \cdots \wedge X_q \wedge R$ does not have divisorial
components   $a_{ij}^l $ is a polynomial.

The identities  (\ref{E:util}) imply that
\begin{eqnarray*}
d \eta &=& \left(L_{X_1} - i_{X_1} d \right) i_{X_2} \cdots i_{X_q}
i_R \Omega
\\&=&\left( i_{[X_1,X_2]} - i_{X_2}L_{X_1} -  i_{X_1} d i_{X_2} \right)i_{X_3}
\cdots i_{X_q} i_R \Omega \, .
\end{eqnarray*}
Using similar manipulations we can proceed by induction on $q$ to
deduce that $d\eta$ is equal to
\[
 (-1)^q  \left(n+1  - q + \sum_i \deg(X_i) \right) i_{X_1}\cdots i_{X_q}  \Omega  +
  \sum_{i} \lambda_i i_{X_1}\cdots \widehat{i_{X_i}}
\cdots i_{X_q}  i_R  \ \Omega \, ,
\]
where the $\lambda_i$  are homogeneous  polynomials of degree
$\deg(X_i) - 1$. If we set
\[ \widetilde X_i =
X_i  + \frac{ \lambda_i R}{(-1)^q  \left(n+1  - q + \sum_i \deg(X_i)
\right) }
\]
then the lemma follows.
\end{proof}

\subsection{The Tangent Sheaf of Foliations}\label{S:sheaf}

Let $\mathcal F$ be a holomorphic foliation  on $\mathbb P^n$, $n\ge
3$, induced by a twisted $q$-form $\omega$. As in the Introduction
the tangent sheaf of $\mathcal F$, denoted by $T \mathcal F$, is the
coherent subsheaf of $T  \mathbb P^n$ generated by the germs of
vector fields annihilating $\omega$.

\medskip

In general $T \mathcal F$ is not locally free. For instance, let
$\mathcal F$ be the codimension one foliation of $\mathbb C^3$
induced by $df$, where $f:\mathbb C^3 \to \mathbb C$ is the function
$f(x,y,z)=x^2+y^2+z^2$. Clearly $T\mathcal F$ is a  locally free
subsheaf of $T\mathbb C^3$ outside the origin of $\mathbb C^3$ since
at these points $f$ is a local submersion. Nevertheless, at the
origin of $\mathbb C^3$, $T \mathcal F$ is not locally free. 

More generally, for codimension one foliations if $T \mathcal F$ is
locally free in a neighborhood $U$ of a point $p$ then the singular
scheme  of $\mathcal F$ on $U$ is defined by the $(n-1)$-minors of a
$n\times (n-1)$ matrix. In particular it is either empty or has
codimension $2$.

\medskip

We say that the tangent sheaf of $\mathcal F$  {\bf splits}  if
\[
   T \mathcal F = \bigoplus_{i=1}^{n-q} \mathcal O_{\mathbb P^n}(e_i)
   ,
\]
for some integers $e_i$. Note that the inclusion of $T\mathcal F$ in
$T \mathbb P^n$ induces sections $X_i \in \mathrm H^0(\mathbb P^n,
T{\mathbb P^n}(-e_i))$ for $i=1\ldots n-q$. It follows from the
Euler sequence that these sections are defined by  homogeneous
vector fields of degree $1-e_{i}\ge 0$ on $\mathbb C^{n+1}$, which
we still denote by $X_i$. The foliation $\mathcal F$ is induced by
the homogeneous $q$-form on $\mathbb C^{n+1}$
\[
  \omega = i_{X_1} \cdots i_{X_{n-q}} i_R \Omega \, .
\]

\subsection{The singular set of $d \omega$ and Kupka Singularities} Another key hypothesis for our results for
codimension one foliations is that $\mathrm{codim} \
\mathrm{sing}(d\omega) \ge 3$. We will now  explain  the {\it
geometrical meaning} of this hypothesis.

If $\omega_0$ is a germ of integrable holomorphic $q$-form on
$(\mathbb C^n,0)$ such that $\omega_0(0)=0$ then  $0$  is called a
Kupka singularity of $\omega_0$ if $d\omega_0(0)\neq 0$. The local
structure of codimension one foliations in a neighborhood of Kupka singularities is
fairly simple: the germ of  foliation is the pull-back of a germ of
foliation on $(\mathbb C^2, 0)$, {\it cf.} \cite[proposition
1.3.1]{Airton} and references therewithin. As a side remark we point
out that the result proved in {\it loc. cit.}  also holds for
integrable $q$-forms $\omega$: if $d\omega(0)\neq0$ then the germ of
codimension $q$ foliation induced by $\omega$ is the pull-back of a
germ of foliation on $(\mathbb C^{q+1}, 0)$.

If $u \in \mathcal O^n_0$ is a unit then $d(u\omega_0) = du
\wedge \omega_0 + u d\omega_0$. Thus the singular set of
$d\omega_0$ is in principle distinct from the singular set of
$d(u\omega_0)$, i.e., it is not an invariant of $\mathcal F_0$,
the foliation induced by $\omega_0$. Although
\[
 \mathcal B(\mathcal F_0)= \mathrm{sing}(\omega_0) \cap \mathrm{sing}(d\omega_0) =
\mathrm{sing}(u\omega_0) \cap \mathrm{sing}(d(u\omega_0)) \, , \] is
a invariant of $\mathcal F_0$ which we will call the non-Kupka
singular set of $\mathcal F_0$.

It is easy to verify that for $\omega$ homogeneous
$1$-form on $\mathbb C^{n+1}$ inducing a foliation $\mathcal F$ of
$\mathbb P^n$ that
\[
\mathcal B(\mathcal F) =  \mathrm{sing}(d\omega) \, .
\]
In other words our hypothesis is on the codimension of the {\it
non-Kupka} singular set of $\mathcal F$.

\section{{Division Lemmata}}

\begin{lemma}\label{L:DRcampos}
Let $X_1, \ldots, X_{n-1}$ be  holomorphic vector fields on $\mathbb
C^{n+1}$ and $\Theta \in \Omega^2(\mathbb C^{n+1})$ be the $2$-form
given by  $ \Theta=i_{X_1} \cdots i_{X_{n-1}} \Omega$. Suppose that
$\mathrm{codim} \ \mathrm{sing}(\Theta) \ge 3$. If $\eta \in
\Omega^2(\mathbb C^{n+1})$ is such that
\[
   \Theta \wedge \eta =0
\]
then there exist holomorphic vector fields $\widetilde X_1, \ldots,
\widetilde X_{n-1}$ such that
\[
   \eta = \sum_{i=1}^{n-1} i_{X_1} \cdots i_{X_{i-1}} i_{\widetilde X_i} i_{X_{i+1}} \cdots
   i_{X_{n-1}}\Omega.
\]
\end{lemma}
\begin{proof}
This follows from  the dual version  of the main result in \cite{Saito};
see also \cite[Proposition (1.1) with $q=3$]{Malgrange}. 
\end{proof}

This Lemma is the case $q=2$ of
Lemma \ref{L:Grass} below (notice that when $q=2$ the Pl\"{u}cker
relations below are $\Theta \wedge \Theta =0$ and the tangent space
of the Grassmannian is given by $\Theta \wedge \eta =0$).

\medskip

\begin{lemma}\label{L:Grass}
Let $X_1, \ldots, X_{n+1-q}$ be  holomorphic vector fields on
$\mathbb C^{n+1}$ and $\Theta \in \Omega^q(\mathbb C^{n+1})$ be the
$q$-form given by $ \Theta=i_{X_1} \cdots i_{X_{n+1-q}} \Omega$.
Suppose that $\mathrm{codim} \ \mathrm{sing}(\Theta) \ge 3$. If
$\eta \in \Omega^q(\mathbb C^{n+1})$ is such that
\[
  {i}_v(\eta) \wedge \Theta +   {i}_v(\Theta) \wedge \eta   =0
\]
for all $v \in \wedge^{q-1} T(\mathbb C^{n+1})$, then there exist
holomorphic vector fields $\widetilde X_1, \ldots, \widetilde
X_{n+1-q}$ such that
\[
   \eta = \sum_{i=1}^{n+1-q} i_{X_1} \cdots i_{X_{i-1}} i_{\widetilde X_i} i_{X_{i+1}} \cdots i_{X_{n+1-q}} \Omega.
\]
\end{lemma}

\begin{proof}
Denote $V = \mathbb C^{n+1}$ and $\mathbf G = \mathrm{Grass}(V, q)$
the Grassmannian of linear subspaces of $V$ of codimension $q$ (i.
e. dimension $n+1-q$). The Pl\"{u}cker embedding
\[
\wp: \mathbf G \to \mathbb P \bigwedge^q(V^*)
\]
gives an isomorphism of $\mathbf G$ with the subvariety of
decomposable $q$-linear forms.
 It is well known (see for example \cite{Griffiths-Harris}) that a $q$-linear form $\theta \in \wedge^q(V^*)$ is decomposable if and only if
\[
({i}_v\theta) \wedge \theta = 0
\]
for every $v \in \wedge^{q-1}V$. It is also known that these equations, 
the well known Pl\"{u}cker quadrics,  generate the ideal of $\wp(\mathbf G)$. Therefore, the
tangent space of $\wp(\mathbf G)$ at a point $\theta$ may be
described as the set of $q$-linear forms $\eta$ such that
\[
 ( {i}_v\eta) \wedge \theta +   ({i}_v \theta) \wedge \eta   = 0
\]
for every $v \in \wedge^{q-1}V$.

Let us consider the standard $q$-multilinear map
\[
\mu: (V^*)^q \to \bigwedge^q(V^*)
\]
defined by $\mu(u_1, \dots, u_q) = u_1 \wedge \dots \wedge u_q$.

If $(V^*)^q_0 = (V^*)^q - \mu^{-1}(0)$ is the open set consisting of
linearly independent vectors then $\mathbf G$ is the quotient of
$(V^*)^q_0$ by the general linear group $\mathrm{GL}(V)$ and the
Pl\"{u}cker embedding is the quotient map (i. e. projectivization) of
$\mu$. Hence the tangent space of $\wp(\mathbf G)$ at a point
$\theta=\mu(u_1, \dots, u_q)$ coincides with the image of the
derivative of $\mu$ at $u=(u_1, \dots, u_q)$, which is given by
\[
d\mu(u).(\widetilde u) = \sum_{i=1}^q u_1 \wedge \cdots \wedge
u_{i-1} \wedge \widetilde u_i \wedge u_{i+1} \wedge \cdots \wedge
u_q
\]
for $\widetilde u=(\widetilde u_1, \dots, \widetilde u_q) \in
(V^*)^q$.

Contraction with $\Omega=dx_0 \wedge \cdots
\wedge dx_n$ induces an isomorphism $\wedge^{n+1-q}(V) \cong
\wedge^q(V^*)$. Therefore  $\wp(\mathbf G)$ is also the projective
image of the multilinear map
\[
\nu: (V)^{n+1-q} \to \bigwedge^q(V^*)
\]
defined by $\nu(v_1, \dots, v_{n+1-q}) = i_{v_1} \cdots
i_{v_{n+1-q}} \Omega$. Hence, the tangent space of $\wp(\mathbf G)$
at a point $\theta$ has a third description as the image of the
derivative of $\nu$
\[
d\nu(v).(\widetilde v) = \sum_{i=1}^{n+1-q} i_{v_1} \cdots
i_{v_{i-1}} i_{\widetilde v_i} i_{v_{i+1}} \cdots i_{v_{n+1-q}}
\Omega
\]
for $v=( v_1, \dots, v_{n+1-q}) \in V^{n+1-q}_0$ and $\widetilde
v=(\widetilde v_1, \dots, \widetilde v_{n+1-q}) \in V^{n+1-q}$.

Returning to the statement that we are proving, let us remark that a
differential $q$-form in $\mathbb C^{n+1}$ may be considered as a
map $\mathbb C^{n+1} \to \wedge^q(\mathbb C^{n+1})^*$. Our
hypothesis about $\Theta$ implies that $\Theta(x)$ is decomposable
for all $x \in \mathbb C^{n+1}$ and then $\Theta$ induces a regular
map
$$\bar \Theta: U \to \mathbf G$$
where $U = \mathbb C^{n+1} - \mathrm{sing}(\Theta)$

The hypothesis on $\eta$ means that $\eta(x)$
belongs to the tangent space to $\mathbf G$ at $\bar \Theta(x)$ for
all $x \in U$. By the last characterization of the tangent space of
$\mathbf G$, there exists an open cover $U = \cup_{\alpha}
U_{\alpha}$ and holomorphic vector fields $\widetilde X_1^{\alpha},
\ldots, \widetilde X_{n+1-q}^{\alpha}$ on $U_{\alpha}$ such that
\[
   \eta|_{U_{\alpha}} =
\sum_{i=1}^{n+1-q} i_{X_1} \cdots i_{X_{i-1}} i_{\widetilde
X_{i}^{\alpha}} i_{X_{i+1}} \cdots i_{X_{n+1-q}} \Omega.
\]

In $U_{\alpha} \cap U_{\beta}$ we have
\[
0 = \eta|_{U_{\alpha}} -\eta|_{U_{\beta}} = \sum_{i=1}^{n+1-q}
i_{X_1} \cdots i_{X_{i-1}} i_{(\widetilde X_{i}^{\alpha} -
\widetilde X_{i}^{\beta})} i_{X_{i+1}} \cdots i_{X_{n+1-q}} \Omega.
\]
Then, for each $j=1, \dots, n+1-q$
\[
i_{X_{j}}(\eta|_{U_{\alpha}} -\eta|_{U_{\beta}}) = \pm
i_{(\widetilde X_{j}^{\alpha} - \widetilde X_{j}^{\beta})} \Theta =
0
\]
Therefore $\widetilde X_{j}^{\alpha} - \widetilde X_{j}^{\beta}$ is
a linear combination of $X_1, \ldots, X_{n+1-q}$ with holomorphic
coefficients. If $X$ denotes the matrix with columns $X_1, \ldots,
X_{n+1-q}$ (and similar notation for $\widetilde X^{\alpha}$), then
there exists a matrix $A^{\alpha \beta}$ with coefficients
holomorphic functions in $U_{\alpha} \cap U_{\beta}$ such that
\[
\widetilde X^{\alpha}  -  \widetilde X^{\beta} = A^{\alpha \beta} X
\]
The collection $\{A^{\alpha \beta}\}_{\alpha \beta}$ is a 1-cocycle
and defines an element of
\[
\mathrm H^1\left(U, \mathcal O_{\mathbb C^{n+1}}^{(n+1)
(n+1)}\right) = \mathrm H^1\big(U, \mathcal O_{\mathbb
C^{n+1}}\big)^{(n+1) (n+1)}  \, .\] The hypothesis $\mathrm{codim} \
\mathrm{sing}(\Theta) \ge 3$ implies that $H^1(U, \mathcal
O_{\mathbb C^{n+1}}) = 0$ (see \cite[pg. 133]{Grauert}). Hence,
after a refinement of the open cover, we may write $A^{\alpha \beta}
= A^{\alpha} - A^{\beta}$ where $A^{\alpha}$ is a holomorphic matrix
in $U_{\alpha}$. Then
\[
\widetilde X^{\alpha} - A^{\alpha} X = \widetilde X^{\beta} -
A^{\beta} X
\]
in $U_{\alpha} \cap U_{\beta}$ and hence the columns of these
matrices define the required  holomorphic vector fields in $U$. To
conclude we apply Hartog's extension Theorem to extend these vector
fields to $\mathbb C^n$.
\end{proof}

\begin{remark}\rm
If the vector fields $X_1, \ldots, X_{n+1-q}$ are homogeneous we can
take the vector fields $\tilde{X}_1, \ldots, \tilde{X}_{n+1-q}$
homogeneous and satisfying $\deg(\tilde{X}_i)=\deg(X_i)$ for all
$i=1\ldots n+1-q$. In fact, if we replace $\tilde{X_i}$ by its
homogeneous component of degree $\deg(X_i)$  then  we still have
that
\[
   \eta = \sum_{i=1}^{n+1-q} i_{X_1} \cdots i_{X_{i-1}} {i_{\tilde X}}_i i_{X_{i+1}} \cdots i_{X_{n+1-q}} \Omega.
\]
\end{remark}

\section{{Foliations with Split Tangent Sheaf}}

In this section we will prove Theorems \ref{T:splits1} and
\ref{T:splits2}.

\begin{proof}[\bf{Proof of Theorem \ref{T:splits1}}] Let $\omega \in
\mathrm{H}^0(\mathbb P^n,\Omega^1_{\mathbb P^n}(d+2))$ be a
saturated  integrable twisted $1$-form on $\mathbb P^n$ and
$\mathcal F$ the induced foliation. If $T\mathcal F$ splits then
there exists a collection $X_1, X_2,\ldots, X_{n-1},R$ of
homogeneous vector fields in involution such that
\[
  \omega = i_{X_1} \cdots i_{X_{n-1}} i_R \Omega \, .
\]
Using lemma \ref{L:prepara} we can also  assume that
\[
    d \omega =  \lambda \cdot i_{X_1} \cdots i_{X_{n-1}}  \Omega \,
    ,
\]
for some $\lambda \in \mathbb C^*$.

 If  $T_{\omega}=
T_{\omega}\mathscr F_1(n,d)$ denotes the Zariski tangent space of
the  scheme $\mathscr F_1(n,d)$ at the point $\omega$ then $\eta \in
T_{\omega}\mathscr F_1(n,d)  $ if, and only if,
\[
   (\omega + \epsilon \eta)\wedge (d\omega + \epsilon d \eta) = 0
   \mod  \epsilon^2  \, .
\]
It follows that  $\eta \in T_{\omega}$ if, and only if, $\omega
\wedge d\eta + \eta \wedge d \omega= 0$.  Note that $\omega \wedge
d\eta + \eta \wedge d \omega$ is annihilated by the radial vector
field. Then  from (\ref{E:dif2}) we conclude that
\[
\omega \wedge d\eta + \eta \wedge d \omega = 0 \, \, \, \,\,
\,\Longleftrightarrow \, \,\, \,\, \, \,\, \,\, d\omega \wedge
d\eta =0 \, .
\]

Thus we can  apply Lemma \ref{L:DRcampos} to $d\eta$ and assure the
existence of a collection $\widetilde X_1, \ldots, \widetilde
X_{n-1}$ of homogeneous vector fields such that
\[
d\eta= \sum_{i=1}^{n-1} i_{X_1} \cdots i_{X_{i-1}} i_{\widetilde
X_i} i_{X_{i+1}} \cdots
   i_{X_{n-1}}\Omega
\]

Since  $i_R\eta=0$ it follows from (\ref{E:dif2})  that
\[
\eta= \sum_{i=1}^{n-1} i_{X_1} \cdots i_{X_{i-1}} i_{\widetilde X_i}
i_{X_{i+1}} \cdots
   i_{X_{n-1}}i_R\Omega \, .
\]

Let $d_i= \deg(X_i)= \deg(\widetilde X_i)$ and denote by $\mathscr
X(d_i)$ the $\mathbb C$-vector space of degree $d_i$ homogeneous
polynomial vector fields on $\mathbb C^{n+1}$. Consider the
alternate multi-linear map
\begin{eqnarray*}
\Psi: \bigoplus_{i=1}^{n-1} \mathscr X(d_i) &\longrightarrow&
\mathrm{H}^0(\mathbb P^n,\Omega^1_{\mathbb P^n}(d+2)) \\
(Y_1,\ldots, Y_{n-1}) &\mapsto& i_{Y_1} \cdots i_{Y_{n-1}} i_R
\Omega \, .
\end{eqnarray*}
The derivative of $\Psi$ at $Y=(Y_1,\ldots, Y_{n-1})$ is
\[
 d \Psi(Y)(Z_1,\ldots, Z_{n-1}) = \sum_{i=1}^{n-1}  i_{Y_1} \cdots i_{Y_{i-1}} i_{{Z_i}}
i_{Y_{i+1}} \cdots
   i_{Y_{n-1}}i_R\Omega \, .
\]

It is now clear  that the every $\eta \in T_{\omega} \mathscr
F_1(n,d)$ is contained in the  image of $d\Psi(X_1, \ldots,
X_{n-1})$. This is sufficient to assure that the image of $\Psi$
contains an open neighborhood of $\omega$ in $\mathscr F_1(n,d)$.
\end{proof}

\medskip

We will see in the proof of Proposition \ref{P:contra} that the
hypothesis on the singular set of $d\omega$ is indeed necessary.

\medskip

The proof of Theorem \ref{T:splits2} is analogous to
the proof of Theorem \ref{T:splits1}; we highlight the unique difference.

\begin{proof}[\bf{Proof of Theorem \ref{T:splits2}}]
If $T_{\omega}= T_{\omega}\mathscr D_q(n,d)$ stands for the Zariski
tangent space of the  scheme $\mathscr D_q(n,d)$ at the point
$\omega$ then  $\eta \in T_{\omega}\mathscr D_q(n,d)$ if, and only
if, ${i}_v(\eta) \wedge \Theta +   {i}_v(\Theta) \wedge \eta   =0$
for all $v \in \wedge^{q-1} T(\mathbb C^{n+1})$. Now we apply the
division lemma \ref{L:Grass} and conclude as in the proof of Theorem
\ref{T:splits1}.
\end{proof}

\section{Two Applications}\label{S:applications}
\subsection{{Linear Pull-backs}}\label{S:Linear}

Our first application is a generalization of a well-known  result by
Camacho and Lins Neto  which says that the pull-back of generic
degree $d$ foliations of $\mathbb P^2$ under generic linear
projections form an irreducible component of $\mathscr F_1(n,d)$ for
every $n\ge 3$, see \cite{CaLN}. More precisely we prove the

\begin{cor}\label{C:PB1}
Let $\mathcal C$ be an irreducible component of $\mathscr F_q(n,d)$.
If the generic element of $\mathcal C$ satisfies the hypothesis of
Theorem \ref{T:splits1} (for $q=1$) or Theorem \ref{T:splits2} (for
$q \ge 2$) then for every integer $m\ge 1$ there exists an
irreducible component of $\mathscr F_q(n + m,d)$ such that the
generic element is the pull-back under a generic linear projection
of a generic element of $\mathcal C$.
\end{cor}
\begin{proof}
Let $\mathcal G  \in \mathscr F_q(n,d)$  be a foliation whose
tangent sheaf splits, i.e.,
\[
 T \mathcal G = \bigoplus_{j=1}^{n-q} \mathcal O_{\mathbb
 P^n}(e_j) \, .
\]
Suppose also that $\mathcal G$ is induced by a $1$-form $\omega$
with $\mathrm{codim} \ \mathrm{sing}(d\omega)\ge 3$ when $q=1$ or
$\mathrm{codim} \ \mathrm{sing}(\omega)\ge3$ when $q\ge 2$. If
$\mathcal F \in \mathscr F_q(n+m,d)$ is the pull-back of $\mathcal
G$ under a generic linear projection $\pi : \mathbb P^{n+m}
\dashrightarrow \mathbb P^n$ then
\[
T \mathcal F= ( \bigoplus_{j=1}^{n-q} \mathcal O_{\mathbb
 P^{n+q}}(e_j) ) \oplus  \mathcal O_{\mathbb
P^{n+q}}(1)^{q}. \] Moreover the codimensions of the singular set of
$\omega$, resp. $d\omega$, and $\pi^* \omega$, resp. $\pi^*
d\omega$, are the same. From Theorems \ref{T:splits1} and
\ref{T:splits2} it is sufficient to prove that every foliation
$\mathcal F' \in \mathscr F_q(n+m,d)$ with $T\mathcal F' = T
\mathcal F$ is the pull-back of a foliation $\mathcal G' \in
\mathscr F_q(n,d)$ under a linear projection.

From the splitting type of $\mathcal F'$ we see that it is induced
by an $i$-form $\omega'$ that may be written as
\[
 \omega'= i_{X_1}\cdots i_{X_{n-q}} i_{Z_1} \cdots i_{Z_q} i_R
 \Omega \, ,
\]
where the $X_j$ are homogeneous vector fields of degree $1-e_j$ and
the $Z_j$ are constant vector fields. In suitable coordinate system
$(z_0, \ldots, z_{n+m})$ of $\mathbb C^{n+m+1}$ we can write $Z_j
=\frac{\partial}{\partial z_{n+j}}$. It follows that the fibers of
the linear projection $\pi'(z_0,\ldots, z_{n+m}) = (z_0, \ldots,
z_n)$ are everywhere tangent to  the leaves of $\mathcal F'$. In
particular the leaves of $\mathcal F'$ are dense open sets of  cones
over the center of projection.  Thus there exists $\mathcal G' \in
\mathscr F_q(n,d)$ such that $\mathcal F'= \pi'^* \mathcal G'$. For
more details the reader may consult \cite[lemma 2.2]{5autores}.
\end{proof}

 As an immediate corollary we obtain irreducible components of $\mathscr
F_q(n,d)$ for arbitrary $q\ge 1$, $n\ge q+2$ and $d\ge 0$.

\begin{example}\label{E:pullback}
The pull-back to $\mathbb P^{n}$ under linear projections of degree
$d$ foliations by curves on $\mathbb P^{q+1}$ fill out irreducible
components of $\mathscr F_q(n,d)$.
\end{example}

Among these irreducible components the only ones that have appeared
before in the literature are the ones with $q=1$ or $d=0$.

\subsection{{Foliations associated to Affine Lie Algebras}}


Foliations induced by representations of the affine Lie algebra in one variable
into the Lie algebra of polynomial vector fields on $\mathbb C^n$, $n\ge
4$,  with homogeneous generators with the usual hypothesis on the
singular set also fill out components of $\mathscr F_q(n,d)$. More
precisely, we prove the

\begin{COR}\label{C:afim}
Let $\mathcal F$ be a codimension $q$ foliation of $\mathbb P^{2+q}$
given by a $q$-form $\omega = i_X i_Y i_R \Omega$ where $\Omega$ is
the euclidean volume form on $\mathbb C^{3+q}$ and $X, Y$ are
homogeneous vector fields of degree $1$ and $e\ge 2$ satisfying the
relation
\[
  [ X , Y ] =  Y   \,.
\]
If $\mathrm{codim} \, \mathrm{sing}(d \omega) \ge 3$ (for  $q=1$) or
$\mathrm{codim} \, \mathrm{sing}( \omega) \ge 3$ (for $q\ge 2$) then
any foliation $\mathcal F'$ sufficiently close to $\mathcal F$ is
induced by a $q$-form $\omega'= i_{X'} i_{Y'} i_R \Omega$ where $X',
Y'$ are homogeneous vector fields of degree $1$ and $e$ satisfying
\[
  [ X' , Y' ] = Y' .
\]
\end{COR}
\begin{proof}
Let $\omega'$ be an integrable $q$-form sufficiently close to
$\omega$. It follows from Theorem \ref{T:splits1} (when $q=1$) and
Theorem \ref{T:splits2} (when $q\ge 2$) that there exist homogeneous
vector fields $X'$ and $Y'$ such that $\omega' = i_{X'}i_{Y'}i_R
\Omega$.

From the integrability of $\omega'$ one deduces that
\[
  [X',Y']= aX' + \lambda Y + bR \, ,
\]
for suitable $\lambda \in \mathbb C$ and $a,b \in S_{e-1}$. Here
 $S_{e-1}$ denotes the space of homogeneous polynomials of
degree $e-1$. Notice that since we can take $X'$ and $Y'$
sufficiently close to $X$ and $Y$ and similarly $\lambda$
sufficiently close to $1$. Moreover, after replacing $X'$ by
$\lambda^{-1}X'$, we can assume that $\lambda=1$.

Since, for arbitrary $\mu \in \mathbb C$,  $$i_{X'}i_{Y'}i_R \Omega
= i_{X'+ \mu R}i_{Y'}i_R \Omega$$ we can suppose that the linear map
\begin{eqnarray*}
T: S_{e-1} &\to& S_{e-1} \, \\
f &\mapsto& X'(f) - f  \,
\end{eqnarray*}
is invertible.

If we set $X''=X'$ and $Y'' = Y' - T^{-1}(a) X' - T^{-1}(b) R$  then
\begin{eqnarray*}
[X'',Y''] &=& aX'+ Y' + bR - X'(T^{-1}(a))X' - X'(T^{-1}(b))R   \\
&=&  (a-X'(T^{-1}(a))X' + Y' + (b-X'(T^{-1}(b))R  \\
&=& Y'' \, .
\end{eqnarray*}
Notice that  $i_{X''} i_{Y''} i_R \Omega = i_{X'} i_{Y'} i_R \Omega$
to conclude the proof of  the corollary.
\end{proof}

When $q=1$ the result below appeared  in a slightly different form
in \cite{CCGL}. For $q\ge 2$ it is new. Using the corollary above
 it is possible to  obtain irreducible
components of $\mathscr F_q(q+2,d)$, for every $q,d$, adapting the
constructions presented in \cite{CCGL}.

\section{{Foliations induced by Group Actions}}\label{S:exemplosLie}

To recall the basic definitions on foliations induced by group
actions, let $\mathfrak g \subset \mathfrak{sl}(n+1, \mathbb C)$ be
a Lie subalgebra of dimension $n-q$. Since $\mathrm H^0(\mathbb P^n,
T{\mathbb P^n}) = \mathfrak{sl}(n+1, \mathbb C)$ is the Lie algebra
of $\mathrm{Aut}(\mathbb P^n) = \mathrm{PSL}(n+1, \mathbb C)$, we
may view $ \bigwedge^{n-q} \mathfrak g $ as a one-dimensional linear
subspace of $ \bigwedge^{n-q} \mathrm H^0(\mathbb P^n, T \mathbb
P^n)$.

By duality we obtain a twisted integrable $q$-form
$$\omega(\mathfrak g) \in
\mathrm H^0(\mathbb P^n, \Omega^{q}_{\mathbb P^n}(n+1)).$$ When
$\omega(\mathfrak g) \neq 0$, the leaves of the foliation $\mathcal
F(\mathfrak g)$ induced by $\omega(\mathfrak g)$ are tangent to the
orbits of $\exp(\mathfrak g)$, the connected (but not necessarily
closed) subgroup of $\mathrm{Aut}(\mathbb P^n)$ with Lie algebra
$\mathfrak g$. Moreover  the  action of $\exp(\mathfrak g)$ on
$\mathbb P^n$  is locally free outside  $\mathrm{sing} (
\omega(\mathfrak g))$.  Notice that
$$\mathrm{deg}(\omega(\mathfrak g)) = n - q = \mathrm{dim}(\mathfrak g)$$

When $\omega(\mathfrak g) \neq 0$ and its  singular set  has no
divisorial components then, clearly,   $T\mathcal F(\mathfrak g) =
\mathfrak g \otimes \mathscr O_{\mathbb P^n}$. In particular, as an
immediate consequence of Theorems \ref{T:splits1} and
\ref{T:splits2} we obtain the

\begin{cor}\label{C:Lie}
Let $\mathfrak g \subset \mathfrak{sl}(n+1, \mathbb C)$ be a Lie
subalgebra of dimension $n-q$ such that $\omega(\mathfrak g)\neq 0$.
If $\omega(\mathfrak g)$ satisfies the hypothesis of Theorem
\ref{T:splits1} (for $q=1$) or Theorem \ref{T:splits2} (for $q \ge
2$) then for every foliation $\mathcal F'$ sufficiently close to
$\mathcal F(\mathfrak g)$ there exists a Lie subalgebra $\mathfrak
g' \subset \mathfrak{sl}(n+1, \mathbb C)$ such that $\mathcal F'=
\mathcal F(\mathfrak g')$.
\end{cor}

Corollary \ref{C:Lie} provides a way to translate  results from the
representation theory of Lie algebras to results about foliations.
Let us consider a simple example.

\begin{example}[Diagonal algebras]\label{E:diagonal}\rm
Let $\mathfrak g$ be a Lie subalgebra of $\mathfrak{sl}(n+1, \mathbb
C)$ of dimension $n - q$ such that  every $g \in
 \mathfrak{g}$ has all eigenvalues with multiplicity one. In
 particular, by an elementary classical result on Lie algebras,
$\mathfrak g$ is diagonal in a suitable system of coordinates.
Moreover from the choice of $\mathfrak g$ this property is clearly
stable: every $\mathfrak g'$ with generators sufficiently close to
the generators of $\mathfrak g$ is also diagonalizable.

If we identify $\mathfrak{sl}(n+1, \mathbb C)$ with the Lie
algebra of linear homogeneous vector fields on $\mathbb C^{n+1}$ with
zero divergence then we can write
\[
\mathfrak{g} \  =  \  < X_1, \ldots, X_{n - q} >
\]
where
\[
X_i = \sum_{j=0}^{n} \lambda_{ij} z_j \frac{\partial}{\partial z_j}
\]
for suitable complex numbers $\lambda_{ij}$. Consequently, we can
associate to $\mathfrak g$ the $\mathcal L$-foliation $\mathcal
F({\mathfrak g}) \in \mathscr F_q(n,n-q)$ induced by the $q$-form
\[
  \omega({\mathfrak{g}}) = i_{X_1}\cdots i_{X_{n-q}} i_R \Omega \, .
\]

When $q=1$ a straightforward computation shows that
\[
d \omega(\mathfrak g) = \pm i_{X_1}\cdots i_{X_{n-1}} \Omega.
\]
In particular the singular set of $d \omega(\mathfrak g)$ is the
union of the sets $\{z_i=z_j=z_k=0\}$ where $i,j,k$ are pairwise
distinct integers in  $\{0, \ldots, n \}$ and hence it has
codimension $3$. For $q \ge 2$ the singular set of $\omega(\mathfrak
g)$ is defined by similar  conditions and has codimension $q$.

We can apply Corollary \ref{C:Lie} to conclude  the existence of
irreducible components of $\mathscr F_q(n,n-q)$ associated to the
diagonal algebras. When $q=1$ these are the well-known logarithmic
components on $\mathbb P^n$ with poles in $n+1$ hyperplanes, cf.
\cite[Corollary 1.19]{CeDe}. Of course we can use Corollary
\ref{C:PB1} to obtain the known result that logarithmic $1$-forms
with poles on  $r\le n+1$ hyperplanes fill out an irreducible
component. \qed
\end{example}

\medskip

Another immediate consequence of Theorems \ref{T:splits1} and
\ref{T:splits2} is the Theorem \ref{T:3}.

\begin{theorem}[Theorem \ref{T:3} of the Introduction]
Let $\mathcal F$ be a codimension $q$ foliation on $\mathbb P^n$
induced by a Lie subalgebra $\mathfrak g \subset \mathfrak{sl}(n+1,
\mathbb C)$ whose corresponding action is locally free outside a
codimension $2$ analytic subset of $\mathbb P^n$. Suppose that
$\mathcal F$ satisfies the hypothesis of Theorem \ref{T:splits1}
(for $q=1$) or Theorem \ref{T:splits2} (for $q \ge 2)$. If $\
\mathrm{H}^1\left(\mathfrak{g}, \ \mathfrak {sl}(n+1, \mathbb C)/
\mathfrak g\right) = 0 \ $ then $\mathcal F$ is rigid, i.e., the
closure of the orbit of such foliation under the automorphism group
of $\mathbb P^n$ is an irreducible component of $\mathscr
F_q(n,n-q)$.
\end{theorem}
\begin{proof}
The main result of \cite{richard} says that a subalgebra $\mathfrak
g \subset \mathfrak{sl}(n+1,\mathbb C)$ for which $\mathrm
H^1(\mathfrak g, \mathfrak{sl}(n+1,\mathbb C)/\mathfrak g) =0$ is
rigid.  The Theorem follows at once combining this result with
Corollary \ref{C:Lie} combined with \cite{richard}.
\end{proof}

\subsection{Foliations induced by Semi-Simple Lie Groups}\label{S:semisimple}
We now turn our attention to foliations associated to semi-simple
Lie subalgebras of $\mathfrak{sl}(n+1,\mathbb C)$.
\begin{cor}\label{C:semisimples}
Let $\mathfrak g \subset \mathfrak{sl}(n+1, \mathbb C)$ be a
semisimple Lie subalgebra of dimension $n-q$, $q \ge 2$, such that
$\omega(\mathfrak g)$ is non-zero. If $\omega(\mathfrak g)$
satisfies the hypothesis of
 Theorem \ref{T:splits2}   then $\mathcal F(\mathfrak g)$ is rigid.
\end{cor}
\begin{proof}
If $\mathfrak g$ is semisimple then for any finite dimensional
$\mathfrak g$-module $M$ one has $\mathrm{H}^1\left(\mathfrak{g}, \
V \right) = 0$, see \cite[Ex. 1.b., Chapter I, paragraph 6, page 102
]{Bourbaki}.  Taking  $M =  \mathfrak {sl}(n+1, \mathbb C)/
\mathfrak g$, the result follows from Theorem \ref{T:3}.
 \end{proof}

\medskip

The  reader may wonder why we have not stated the Corollary
\ref{C:semisimples}  for codimension one foliations since the proof
works  also in this case. The reason is very simple: the hypothesis
of Theorem \ref{T:splits1} are never satisfied. More precisely we
have the

\begin{prop}\label{P:contra}
If $\mathfrak g \subset \mathfrak{sl}(n+1, \mathbb C)$ is a
semisimple  Lie subalgebra of dimension $n-1$ such that
$\omega(\mathfrak g)\neq 0$ then $\mathrm{codim} \, \, \mathrm{sing}
(d(\omega(\mathfrak g)) \le 2$.
\end{prop}
\begin{proof}
It follows from \cite[Thm 1.22]{CeDe} that $\mathcal F(\mathfrak g)$
admits a rational integral $F:\mathbb P^n \dashrightarrow \mathbb
P^1$. From Stein's factorization Theorem we can assume that the
generic fiber of $F$ is irreducible and that such $F$ is unique up
to right composition with elements in $\mathrm{Aut}(\mathbb P^1)$.
The argument used to prove the above mentioned result shows that
every fiber of $F$ is irreducible. The point is that for every $X
\in \mathfrak g$ and every $\mathcal F(\mathfrak g)$-invariant
hypersurface $\{P=0\}$
\[
 X(P)=0 \, .
\]
Otherwise it would exist a non-trivial morphism of Lie algebras
$\mu_X: \mathfrak g \to \mathbb C$, cf. loc. cit. for more details.
If one of the fibers of $F$ admits a prime decomposition of the form
$P_1^{n_1} \cdots P_l^{n_l}$ with $l\ge 2$ then
\[
 X\left( \frac{P_1^{\deg(P_2)}}{P_2^{\deg(P_1)}} \right) =0 \, ,
\]
contradicting the unicity of $F$.

We recall that a classical Theorem of Halphen says that a pencil on
$\mathbb P^n$ with irreducible generic element has at most two
multiple elements, cf. \cite{coloquio}.  Since every fiber of
$F$ is irreducible it follows that $F$ has at most two non-reduced
fibers.  In particular, after composing with an automorphism of
$\mathbb P^1$, we can assume that $F$ is of the form
\[
   F =\frac{G^{\deg(H)}}{H^{\deg(G)}}
\]
and every fiber of $F$ distinct from $F^{-1}(0)$ and
$F^{-1}(\infty)$ is reduced and irreducible. It follows that
$\omega(\mathfrak g)$ is a complex multiple of the $1$-form
\[
   \deg(H) H dG - \deg(G) G dH \ .
\]
Notice that $\deg(H) + \deg(G) -2 = \deg(\mathcal F(\mathfrak g))=
n-1$. If we take $H'$ and $G'$ homogeneous polynomials, arbitrarily
close to $H$ and $G$ respectively, that cut out smooth hypersurfaces
intersecting transversely on $\mathbb P^n$ then it follows from
\cite{singularidades} that the  $1$-forms $\deg(H') H' dG' -
\deg(G') G' dH'$ do have  isolated singularities.  In particular the
tangent sheaf of the induced foliations is not locally free. It
follows from Theorem \ref{T:splits1} that $\mathrm{codim} \,
\mathrm{sing} (d(\omega(\mathfrak g)) \le 2$.
\end{proof}

\medskip
The next example shows that in the case of higher codimension we do
have foliations satisfying the hypothesis of Corollary
\ref{C:semisimples}.

\begin{example}[Exceptional component of $\mathscr
F_q(q+3,3)$]\label{E:PSL} If $q\ge 2$ then there exists an
irreducible component of $\mathscr F_q(q+3,3)$ such that the generic
element is conjugate by an automorphism of $\mathbb P^{q+3}$ to the
foliation induced by the natural action of $\mathrm{PSL}(2,\mathbb
C)$ on $\mathrm{Sym}^{q+3} \mathbb P^1\cong \mathbb P^{q+3}.$
\end{example}
\begin{proof}
For each natural number $r=q+3$ let us consider the action of
$\mathrm{PSL}(2,\mathbb C)$ on the projective space $\mathbb P
S^r(\mathbb C^2)^* = \mathbb P^r$ of binary forms of degree $r$. Let
$\mathcal F_r$ be the three-dimensional foliation on $\mathbb P^r$
induced by this action.

A positive divisor $D$ on $\mathbb P^1$ has finite stabilizer if,
and only if, its support contains at least three points. Hence, the
singular set of $\mathcal F_r$ is the union of the two-dimensional
varieties $S_m= \{m p  + (r-m) q, \ p, q \in \mathbb P^1\}$ for $0
\le m \le r/2$. Hence, it follows from Corollary \ref{C:semisimples}
that for $r \ge 5$, equivalently $q\ge 2$,  the foliation $\mathcal
F_r$ is rigid.
\end{proof}

The example above is in fact a particular case of the more general

\begin{example}\label{E:semisimples} Let $G = SL(n, \mathbb C)$
and consider the natural action on the $m$-th symmetric power $V =
S^m(\mathbb C^n)$. More generally, let $G$ be a classical simple Lie
group of dimension $d$ and let $V$ be a finite direct sum of irreducible
representations, for instance, symmetric or alternating powers
of the standard representation. In most of these cases the
hypothesis above on stabilizers is satisfied and hence one obtains
irreducible components of $\mathscr F_{n-d}(n, d)$ corresponding to
these rigid foliations.
\end{example}
\begin{proof}
The hypothesis on the stabilizers implies that $\omega(\mathfrak g)$
satisfies the hypothesis of Theorem \ref{T:splits2}. Hence the
statement follows from Corollary \ref{C:semisimples}.
\end{proof}

\subsection{An Infinite Family of Rigid
Foliations}\label{S:infinito}

 Since the core of Theorem \ref{T:infinito} proof
 consists in stablishing  the vanishing of a certain
cohomology group   we will briefly recall the definition of Lie
algebra cohomology thinking on reader's ease.

If $\mathfrak g$ is a Lie algebra and $M$ is  a $\mathfrak g$-module
then the cohomology groups  $\mathrm H^*(\mathfrak g, M)$ are
defined as the cohomology of the complex $({C}^{*}(\mathfrak g,
M),d)$ in which the $n$-cochains $f\in {C}^{n}(\mathfrak g,M)$ are
multilinear antisymmetric maps
\[ f: \underbrace{\mathfrak{g}\times \ldots \times \mathfrak{g}}_{n\ \text{times}} \longrightarrow M\]
 and the coboundary is given by
\begin{multline*}
 df(v_0,\dots,v_{n})=
\sum_{i=0}^n (-1)^i [v_i,f(v_0,\dots,\hat{v_i},\dots,v_{n})]+\\
+\sum_{i<j} (-1)^{i+j}
f([v_i,v_j],v_0,\dots,\hat{v_i},\dots,\hat{v_j},\dots,v_{n}).
\end{multline*}

\medskip

\begin{proof}[\bf{Proof of Theorem \ref{T:infinito}}]  We are considering the foliation $\mathcal F \in
\mathscr F_1(n,n-1)$, $n\ge3$, induced by the subalgebra $\mathfrak
g \subset \mathfrak{sl}(n+1,\mathbb C)$ generated by
\begin{eqnarray*}
X &=& \sum_{i=0}^n (n-2i) z_i \frac{\partial z}{\partial z_i} \, ,
\\
Y_k &=& \sum_{i=0}^{n-k} z_{i+k} \frac{\partial z}{\partial z_i}\,
,\quad k=1,\ldots, n-2  \, .
\end{eqnarray*}
Notice that
\[
   [X,Y_i] = -2i Y_i \, , \text{ when } i=1,\ldots, n-2, \quad \text{ and }  \quad
   [Y_i,Y_j] = 0 \, \text{ for arbitrary } i, j.
\]
Using (\ref{E:util}) of \S\ref{S:calculus} to compute  $d\omega$
like in lemma \ref{L:prepara} we verify that
\[
d\omega = (-1)^{n-2}
 i_{Z} i_{Y_1} \cdots i_{Y_{n-2}}\Omega \,
 ,
\]
where  $\Omega=dz_0 \wedge \cdots \wedge dz_n$ and  $Z = (n+1)X -
(n-1)(n-2)R$.

The singular locus of $d\omega$ is thus defined by the vanishing of
the $(n-1)\times(n-1)$ minors of the $(n-1)\times(n+1)$ matrix
\[
\left(
  \begin{array}{ccccccc}
    \lambda_0 z_0 & \lambda_1 z_1 &   \cdots & \lambda_{n-1} z_{n-3} & \lambda_{n-2} z_{n-2} & \lambda_{n-1} z_{n-1} &  \lambda_{n}z_n \\
    z_1 & z_2 &  \cdots & z_{n-2} & z_{n-1}  & z_{n} & 0 \\
    z_2 & z_3 &   \cdots & z_{n-1} & z_n & 0 & 0 \\
    \vdots & \vdots & \ldots & \vdots & \vdots & \vdots &\vdots \\
    z_{n-2} &  z_{n-1} &  \cdots  & 0 &  0 & 0 & 0
     \\
  \end{array}
\right) ,
\]
where $\lambda_i = (n+1) (n-2i) - (n-1)(n-2).$

Observe that $\lambda_n$, $\lambda_{n-1}$ and $\lambda_{n-2}$ are
all different from zero when $n\ge 3$. Thus omitting the first two
columns of the matrix above we see that $z_{n}^{n-1}$ appears in the
ideal generated by the $(n-1)\times (n-1)$ minors. Therefore, set
theoretically, $\mathrm{sing}(d\omega) \subset \{ z_n=0\}$. If we
set $z_n=0$  and omit the first and the last column we see that
$\mathrm{sing}(d\omega) \subset \{ z_n=0\} \cap \{z_{n-1}=0\}$.
Analogously after omitting the last two  columns and setting
$z_n=z_{n-1}=0$ we conclude that $\mathrm{sing}(d\omega) \subset \{
z_n=0\} \cap \{z_{n-1}=0\} \cap \{z_{n-2} = 0\}$ for every  $n\ge
3$. In particular  $\mathrm{codim} \, \mathrm{sing}(d\omega)\ge 3$
when $n\ge3$.

\medskip

We have just verified that to prove the rigidity of $\mathcal F$ we
can apply Theorem \ref{T:3} once we know that $\mathrm H^1(\mathfrak
g,   M )=0$ where $M$ is the $\mathfrak g$-module
$\mathfrak{sl}(n+1, \mathbb C) / \mathfrak g$. The remaining part of
the proof will be devoted to establish the vanishing of this
cohomology group.

\medskip

Observe that $\mathrm{ad}(X):M\to M$ is semi-simple and that
\[
M = \bigoplus_{i=-n}^{n} M_{2i}
\]
where $M_{i}$ is the $\mathrm{ad}(X)$-eigenspace corresponding to
the eigenvalue $i$.

Let $f \in C^1(\mathfrak g, M)$ be a cocycle, i.e., $df =0$. We can
assume that $f(X) \in M_0$. Indeed, if $f(X) \notin M_0$ then we
have  just to  replace $f$ by $f-dv$ for a suitable $v \in
C^0(\mathfrak g,M)$.

Remark that $df(X,Y_1)=0 $ implies
\begin{equation}\label{L:col}
[X,f(Y_1)]= - 2f(Y_1)  + [Y_1,f(X)] \, . \end{equation}
 If we write
\[
f(Y_1) = \sum_i v_i \, , \quad \text{ where } v_i \in M_i,
\]
then (\ref{L:col}) can be rewritten as
\[
\left\lbrace \begin{array}{lclcl} \mathrm{ad}(X)(v_i) + 2v_i &=& 0
 &\text{when}& i \neq -2, \\
\mathrm{ad}(X)(v_{-2}) + 2 v_{-2} &=&   \mathrm{ad}(Y_1)(f(X)) \, .
\end{array}
\right.
\]
It is now clear  that
\[
[Y_1,f(X)]=0 \quad \text{ and } \quad f(Y_1) \in M_{-2} \, .
\]

Since $\mathrm{ad}(Y_1)$ maps $M_0$ isomorphically to $M_{-2}$ we
can assume that $f(X)=f(Y_1)=0$. As before we have just to replace
$f$ by $f-dv$ for a suitable $v\in C^0(\mathfrak g,M)$.

Since $f(X)=0$ the identity $d f(X,Y_i)=0 $  holds and consequently
$f(Y_i) \in M_{-2i}$ for all $i \in \{2, \ldots, n-2\}.$

\medskip

Let now $F:\mathfrak g \to \mathfrak{sl}(n+1,\mathbb C)$ be a a
linear map lifting $f:\mathfrak g \to M$ such that the image of $F$
is contained in a vector subspace $\overline M$ of
$\mathfrak{sl}(n+1,\mathbb C)$ satisfying
\[  ad(X)(\overline M)\subset \overline M \quad \text{ and } \quad \mathfrak{sl}(n+1,\mathbb
C)= \overline M \oplus \mathfrak{g} \, .
\]
The existence of  $\overline M$  follows at once from the fact that
$ad(X):\mathfrak{sl}(n+1,\mathbb C) \to \mathfrak{sl}(n+1,\mathbb
C)$ is semi-simple.

If we set $F(Y_k)=B_k $, $k=2,\ldots, n-2$, then $f(Y_k) \in
M_{-2k}$ implies that
\[
  B_k = \sum_{i=0}^{n-k} b_i^{(k)} z_{i+k} \frac{\partial}{\partial
  z _i}\,, \quad \forall k\in \{2,\ldots, n-2\}.
\]
We claim that $B_{n-2}=0$. In fact $df(Y_1,Y_{n-2}) =0$ implies that
\[
  [Y_1, B_{n-2}] = \sum_{i}^1 \left(b_i^{(n-2)} - b_{i+1}^{(n-2)}\right) z_{i+n-1} \frac{\partial}{\partial z_i} = 0 \quad \mod \mathfrak g \, .
\]
Thus $B_{n-2}$ must be a complex multiple of $Y_{n-2}$. Since
$Y_{n-2} \notin \overline M$ the claim follows.

\medskip

To conclude the proof of the Theorem we will show that $B_k=0$ for
every $k \in \{2,\ldots, n-2\}$. Clearly it suffices to settle that
\begin{enumerate}
\item[(a)] If  $B_{n-k}=0$ then $B_k = 0$;
\item[(b)] If $B_k =0$ then $B_{n-(k+1)} = 0$.
\end{enumerate}

To prove that (a) holds  first observe that  $df(Y_1,Y_k) = 0 $
implies that there exists $\lambda_k \in \mathbb C$ such that
$[Y_1,B_k]=\lambda_k Y_{k+1}$. In more explicit terms
\[
  \sum_{i=0}^{n-(k+1)} \left(b_i^{(k)} - b_{i+1}^{(k)}\right) z_{i+k+1} \frac{
  \partial}{\partial z_i} = \lambda_k Y_{k+1} \, .
\]
Thus the sequence $\{ b_i^{(k)} \}_{i=0}^{n-k}$ is an arithmetic
progression with step $-\lambda_k$.

Since $B_{n-k}=0$ and $df(Y_k,Y_{n-k}) = 0 $ then
$[Y_{n-k},B_k]=b_0^{(k)} - b_{n-k}^{(k)}z_{n}\frac{\partial
}{\partial z_0} = 0$. Therefore $\lambda_k =0$ and consequently
$B_k$ is a complex multiple of $Y_k$. Since $Y_k \notin \overline M$
we conclude that  $B_k=0$. Assertion (a) follows.

To prove (b) we will proceed similarly. On the one hand
$df(Y_1,Y_{n-k+1})=0$ implies  that the sequence $\{ b_i^{(n-(k+1))}
\}_{i=0}^{k}$ is an arithmetic progression with step
$-\lambda_{k+1}$. On the other hand $B_k=0$ implies that
$df(Y_k,Y_{n-(k+1}) = [Y_k,B_{n-(k+1)}]=0$. Since
\[
[Y_k,B_{n-(k+1)}] = \sum_{i=0}^{1}\left( b_{i}^{(n-(k+1))} -
b_{i+k}^{(n-(k+1))}\right) z_{i+n-1} \frac{\partial}{\partial z_i}
\,
\]
we obtain that $\lambda_{n-(k+1)} =0$ and that $B_{n-(k+1)}$ is a
complex multiple of $Y_{n-(k+1)}$. Since $Y_{n-(k+1)} \notin
\overline M$ the assertion (b) follows and so does the Theorem.
\end{proof}

\subsection{Two other Rigid Foliations}\label{S:AFF}

The  reader will notice that with minor modifications the proof of
Theorem \ref{T:infinito} also shows that the Lie subalgebras
$\mathfrak g(n,r) \subset \mathfrak{sl}(n+1,\mathbb C)$ generated by
\begin{eqnarray*}
X &=& \sum_{i=0}^n (n-2i) z_i \frac{\partial z}{\partial z_i} \, ,
\\
Y_k &=& \sum_{i=0}^{n-k} z_{i+k} \frac{\partial z}{\partial z_i}\,
,\quad k=1\ldots \mathbf{n-r}  \, ,
\end{eqnarray*}
satisfy $\mathrm H^1(\mathfrak g(n,r), \mathfrak{sl}(n+1, \mathbb C)
/ \mathfrak g(n,r)) = 0$ for every $r \in \{2,\ldots, n-1\}$. The
rigidity of the corresponding foliations will follows from Theorem
\ref{T:3} once we verify that the singular set has codimension at
least $3$.  When the algebra above has dimension two ($r=n-1$) we
are in a particularly interesting situation described in the example
below.

\begin{example}[Exceptional component of $\mathscr
F_q(q+2,2)$]\label{E:AFF} If $q\ge 2$ then there exists an
irreducible component of $\mathscr F_q(q+2,2)$ such that the generic
element is conjugate by an automorphism of $\mathbb P^{q+2}$ to the
foliation induced by the natural action of  $\mathrm{Aff}(\mathbb
C)$ on $\mathrm{Sym}^{q+2} \mathbb P^1\cong \mathbb P^{q+2}.$
\end{example}
\begin{proof}
Let $q \ge 2$ and  $\mathcal F_q$ be the foliation of
$\mathrm{Sym}^{q+2} \mathbb P^1 \cong \mathbb P(\mathbb
C[x,y]_{q+2}) \cong \mathbb P^{q+2}$ induced by the natural action
of the following subgroup of $\mathrm{PSL(2,\mathbb C)} \cong
\mathrm{Aut}(\mathbb P^1)$:
\[
\left\{ \left(
          \begin{array}{cc}
            a & b \\
            0 & a^{-1} \\
          \end{array}
        \right) \Big\vert \, \, a \in \mathbb C^*, b \in \mathbb C \right\}.
\]

A positive divisor $D$ on $\mathbb P^1$ has finite stabilizer if,
and only if, its support contains at least two points of $\mathbb
P^1 - \{\infty\} = \{ [x:y] \in \mathbb P^1 \vert y\neq 0 \}$.
Therefore the generic orbit has dimension two and the singular set
of $\mathcal F_q$ is the union of the one-dimensional varieties
\[
C_m = \{ (q-m) \infty + m p \, ; \, p \in \mathbb P^1 \} \, \, \,
\text{ for } \, \, \, 0\le m \le q \, .
\]
Observe that $C_q$ is the Veronese curve of degree of $q$ in
$\mathbb P^q$ and $C_m$ is a Veronese curve of degree $m$ in the
osculating $\mathbb P^m$ to $C_q$ at the point $\infty$. In
particular the singular set of $\mathcal F_q$ has codimension at
least $3$.

To verify the rigidity of $\mathcal F_r$ we first  give explicit
expressions for the vector fields on $\mathbb P^r$ inducing
$\mathcal F_r$. Since
\[
\begin{array}{lclcl}
  ((1+\epsilon)x)^i ( (1+\epsilon)^{-1}y)^j  & = & x^i y^j + \epsilon (i-j) x^iy^j & \mod  & \epsilon^2 \\
  (x+\epsilon)^i y^j  & = & x^iy^j + \epsilon i x^{i-1}y^{j+1} & \mod  &
  \epsilon^2
\end{array}
\]
it follows that on the basis $< z_i = x^{q-i}y^i
>$ of  $\mathbb C_q[x,y]$ the tangent sheaf of $\mathcal F_q$ is generated by the vector
fields
\[
X=  \sum_{i=0}^{q+2} (q+2-2i) z_i \frac{\partial}{\partial z_i}
\quad \text{and} \quad  Y= \sum_{i=0}^{q+1} (k-i) z_{i+1}
\frac{\partial}{\partial z_i} .
\]
After a change of coordinates of the form $(z_0,\ldots, z_{q+2})
\mapsto (\lambda_0 z_0,\ldots, \lambda_{q+2} z_{q+2})$, where
$(\lambda_0,\ldots, \lambda_{q+2}) \in \mathbb C^{q+3}$, we can
assume that
\[
X=  \sum_{i=0}^{q+2}  (k-2i) z_i \frac{\partial}{\partial z_i} \quad
\text{and} \quad  Y= \sum_{i=0}^{q+1} z_{i+1}
\frac{\partial}{\partial z_i} .
\]
Thus the corresponding algebra is isomorphic  to $\mathfrak
g(q+2,q+1)$ and the rigidity follows from Theorem \ref{T:3}.
\end{proof}

\medskip

We can thus interpret the foliations obtained in Theorem
\ref{T:infinito} as {\it extensions} of the foliations on $\mathbb
P^n = \mathrm{Sym}^n \mathbb P^1$ induced by the natural action of
$\mathrm{Aff}(\mathbb C)$. Here, by an extension of a foliation
$\mathcal F$ we mean a foliation $\mathcal G$ such that 
$T \mathcal F \subset T \mathcal G$.  Since
the codimension one foliations in question are rigid it is therefore
natural to wonder if these extensions are unique. Below we present
some examples in dimensions $6$ and $7$ showing that this is not the
case. As we will see they also correspond to rigid foliations.

\medskip

To construct the examples we will take $X,Y_1,Y_{n-2} \in
\mathfrak{sl}(n+1,\mathbb C)$ as in Theorem \ref{T:infinito} and
will look for $Y_2,Y_3, \ldots, Y_{n-3}\in \mathfrak{sl}(n+1,\mathbb
C)$ such that for every $k \in \{2,\ldots, n-3\}$ the following
relations holds:
\begin{equation} \label{L:condicoes}
\mathrm{(a)} \,   ad(X)(Y_k) = -2kY_k  \quad \text{and} \quad
\mathrm{(b)} \,   ad(Y_1)(Y_k) = -Y_{k+1}.
\end{equation}

From (\ref{L:condicoes}.a) it follows that $Y_k$ must be of the form
\[
Y_k = \sum_{i=0}^{n-k} b_i^{(k)} z_{i+k} \frac{\partial}{\partial
z_i} \, \quad \text{ for some } b_i^{(k)}\in \mathbb C.
\]
The relations (\ref{L:condicoes}.b) imply that $ b_{i+1}^{(k)} =
b_i^{(k+1)}  + b_i^{(k)}$ for every $i\in \{0,\ldots,n-k-1\}$ and
$k\in \{ 2,\ldots,n-3\}$. It is then an amusing exercise to deduce
that
\begin{equation}\label{E:milagre}
b_i^{(n-k)} = \sum_{l=0}^{k-2} \binom{i}{l} b_0^{(n-k+l)} \quad
\forall k \in \{ 3, \ldots, n-2 \} \, .
\end{equation}

The equations quoted in (\ref{L:condicoes}) together with Jacobi's
relation implies that $[Y_i,Y_j]$ is an eigenvector of $ad(X)$ with
eigenvalue $-2(i+j)$. Thus if the vector subspace of
$\mathfrak{sl}(n+1,\mathbb C)$  spanned by $X,Y_1,\ldots, Y_{n-2}$
is a Lie subalgebra then $Y_i$ also satisfies the relations
\[
   [Y_{n-k}, Y_k ]  = 0 \quad \text{and} \quad     [Y_{n-k-1}, Y_k ]
   = 0\, .
\]
Now notice that
\begin{equation}\label{E:33}
\left\lbrace \begin{array}{lclcl} \,[Y_{n-k},Y_k] =0 &\implies&
b_{i+k}^{(n-k)} b_i^{(k)} - b_i^{(n-k)}
b_{i+n-k}^{(k)} &=& 0 ,\\
\,[Y_{n-k},Y_{k-1}] =0 &\implies& b_{i+k-1}^{(n-k)} b_i^{(k-1)} -
b_i^{(n-k)} b_{i+n-k}^{(k-1)} &=& 0\, .
\end{array}\right.
\end{equation}

The solutions of the system defined through the equations
(\ref{E:milagre}, \ref{E:33}) are completely determined by the
values of $b_0^{(k)}$ where $k$ ranges from $2$ to $n-3$. For $n \in
\{ 5,6,7,8 \}$ we carried out  in detail the calculations. We
summarize below the results:
\begin{itemize}
\item[$\mathbf{n=5}.$] There are no solutions.
\item[$\mathbf{n=6}.$] \noindent There is only one solution. Namely
$
 \left( b_0^{(2)},  b_0^{(3)} \right) =\left( \frac {9}{8},-\frac{3}{2}\right).
$ This solution corresponds indeed to a Lie algebra since
$[Y_2,Y_3]=0$ thanks to (\ref{E:33}).  We will denote the
corresponding subalgebra by $\mathfrak g_6$.
\item[$\mathbf{n=7}.$] There is only one solution. Namely
\[
\left( b_0^{(2)},  b_0^{(3)}, b_0^{(4)} \right)  = \left(
\frac{\sqrt {3}}{2}, 1-\sqrt {3}, \frac{-3+\sqrt {3}}{2} \right)
\]
The only bracket whose vanishing is not imposed by (\ref{E:33}) is
$[Y_2,Y_3]$. It turns out that $$[Y_2,Y_3]= \frac{5}{2}Y_5.$$
 We will denote the corresponding
subalgebra by $\mathfrak g_7$.
\item[$\mathbf{n=8}.$] Here we have  two  possibilities for $\left( b_0^{(2)},  b_0^{(3)}, b_0^{(4)} , b_0^{(5)} \right)$.  Namely
\[
\begin{array}{l}
\left( \frac{45 + 15\sqrt{265}}{256} , \frac{- 15 +5
\sqrt{265}}{64}, \frac{ 35-\sqrt{265}}{32}, -\frac{3}{2} \right)
\,\,\, \text{and} \,\,\,
  (0,0,-1,0).
\end{array}
\]
In  both cases  we have that $[Y_2,Y_3]$ is not a complex multiple
of $Y_5$. Thus they do not correspond to Lie subalgebras.
\end{itemize}

\medskip

\begin{prop}\label{P:67}
The foliations  $\mathcal F_k=\mathcal F(\mathfrak g_k) \in \mathscr
F_1(k,k-1)$, $k=6,7$, are rigid.
\end{prop}
\begin{proof}
Arguing as in the proof of Theorem \ref{T:infinito} we can verify
 that $\mathrm{codim} \, \mathrm{sing} ( d \omega(\mathfrak
g_k) ) \ge 3$ for $k=6,7$. Instead of computing the relevant
cohomology groups we will prove the rigidity of $\mathcal F_6$ and
$\mathcal F_7$ by a more elementary argument.

 Corollary \ref{C:Lie} implies
that every foliation $\mathcal F=[\omega]$ sufficiently close to
$\mathcal F_k$ is induced by a $\mathfrak g \subset
\mathfrak{sl}(k+1,\mathbb C)$. Moreover we can also assume that
$\mathrm{codim} \, \mathrm{sing} ( d \omega ) \ge 3$.

Notice that  $\mathfrak h_k = [\mathfrak g_k, \mathfrak g_k]$ has
codimension one in $\mathfrak g_k$. By semi-continuity it follows
that $\mathfrak h = [\mathfrak g, \mathfrak g]$ has either
codimension one or zero in $\mathfrak g$. Since $\mathrm{codim} \,
\mathrm{sing} ( d \omega ) \ge 3$ it follows from Proposition
\ref{P:contra} that $\mathfrak h$ has indeed codimension one in
$\mathfrak g$.

Let $X' \in \mathfrak g - \mathfrak h$ be sufficiently close to $X
\in \mathfrak g_k$. Since $ad(X):\mathfrak h_k \to \mathfrak h_k$ is
semi-simple we distinct eigenvalues the same holds for $ad(X'):
\mathfrak h \to \mathfrak h$. Thus there exists $Z_1 \in \mathfrak
h$ such that $[X',Z_1]$ is a multiple of $Z_1$ and $Z_1$ is a
deformation of $Y_1$. It follows from Example \ref{E:AFF} that after
a change of coordinates we suppose that $X'=X$ and $Z_1=Y_1$. In
particular the eigenvalues of $ad(X'):\mathfrak h \to \mathfrak h$
are integers and by continuity they are equal to  $ -2,-4,-6,-8$ for
$k=6$ and $ -2,-4,-6,-8,-10$  for $k=7$. Let $Z_1,Z_2,\ldots,
Z_{k-2} \in \mathfrak h_k$ be the corresponding eigenvectors.

Now from Jacobi's identity we deduce that
\[
[X,[Z_i,Z_j]] = -2(i+j) Z_{i+j} \implies [Z_1,Z_j] = \lambda_j
Z_{j+1}, \, \, \, j=2,\ldots, k-2\, .
\]
Consequently after replacing $Y_j$ by a complex multiple for
$j=2,\ldots, k-2$  we can assume that $Z_{k-2}=Y_{k-2}$ and that
$[Z_1,Z_j]=Z_{j+1}$ for $j=2,\ldots, k-3$. The proposition follows
from the calculations made before its statement.
\end{proof}

If $\mathfrak g \subset \mathfrak{sl}(n+1, \mathbb C)$ is a rigid
Lie subalgebra then, in general, we cannot guarantee that $
   \mathrm{H}^1\left(\mathfrak{g}, \
{\mathfrak {sl}(n+1, \mathbb C)}/{\mathfrak g}\right) = 0 \, . $ The
point is that it may happen that the variety of Lie subalgebras is
non-reduced at $\mathfrak g$. This does not happen in the examples
that we studied and we are not aware of any concrete example.
Although R. Carles constructed several examples of rigid Lie
algebras (of dimension at least $9$) where the variety of Lie
algebras is non-reduced, see for instance \cite{Carles1} and
references there within. Thus it its natural to expect that the
$\mathscr F_q(n,d)$ are non-reduced in general. It would be
interesting to construct examples of irrreducible components which
are everywhere non-reduced.

Another intriguing fact is that up to now all the known irreducible
components of $\mathscr F_q(n,d)$ are unirational varieties. It
would be interesting to know if this is a general fact or if it is
just a testimony of our limited knowledge  about  the irreducible
components of the space of holomorphic foliations on projective
spaces.

\vspace {5 mm}

\small
\begin{flushright}
\begin{tabular}{ l  l l l    l }
  {\bf Fernando Cukierman} &&&& {\bf Jorge Vit\'{o}rio Pereira}  \\
  email: fcukier@dm.uba.ar &&&& email: jvp@impa.br \\
  Univ. de Buenos Aires &&&& IMPA \\
  Dto. Matem\'{a}tica, FCEN &&&& Estrada Dona Castorina,110 \\
  Ciudad Universitaria &&&&   22460-320  Jardim Bot\^{a}nico \\
  (1428) Buenos Aires &&&& Rio de Janeiro \\
  Argentina &&&& Brasil
\end{tabular}
\end{flushright}


\begin{thebibliography}{99}


\bibitem{Bourbaki}
N. Bourbaki, \emph{Lie groups and Lie algebras, chapters 1-3}
Springer-Verlag.

\bibitem{Omegar}
O.  Calvo-Andrade, \emph{Irreducible components of the space of
holomorphic foliations.} Math. Ann. \textbf{299} (1994), no. 4,
751--767.


\bibitem{CCGL}
O.  Calvo-Andrade, D. Cerveau, L. Giraldo and A.  Lins Neto,
\emph{Irreducible components of the space of foliations associated
to the affine Lie algebra.}  Ergodic Theory Dynam. Systems
\textbf{24} (2004), no. 4, 987--1014.


\bibitem{CaLN}
C. Camacho and A.  Lins Neto, \emph{The topology of integrable
differential forms near a singularity.} Inst. Hautes \'{E}tudes Sci.
Publ. Math. No. \textbf{55}, (1982), 5--35.


\bibitem{Carles1}
R. Carles, \emph{D\'{e}formations dans les sch\'{e}mas d\'{e}finis par les
identit\'{e}s de Jacobi.}  Ann. Math. Blaise Pascal \textbf{3} (1996),
no. 2, 33--62.



\bibitem{CeDe}
D. Cerveau and J. Deserti, \emph{Feuilletages et actions de
groupes sur les espaces projectifs}, Pr\'{e}publication IRMAR, (2004).


\bibitem{CL}
D. Cerveau and A. Lins Neto, \emph{Irreducible components of the
space of holomorphic foliations of degree two in $CP(n)$, $n\geq
3$}.  Ann. of Math., \textbf{143} (1996), 577--612.

\bibitem{CLE}
D. Cerveau, A. Lins Neto and S. J. Edixhoven, \emph{Pull-back
components of the space of holomorphic foliations on ${\mathbb
C}{\mathbb P}(n)$, $n\geq 3$.}  J. Algebraic Geom.  \textbf{10}
(2001), no. 4, 695--711.

\bibitem{5autores}
D. Cerveau, A. Lins Neto, F. Loray, J. V. Pereira  and F. Touzet,
\emph{Algebraic Reduction Theorem for complex codimension one
singular foliations}, Comment. Mat. Helv., \textbf{81} (2006),
157-169.

\bibitem{singularidades}
F. Cukierman, M. Soares and I. Vainsencher, \emph{Singularities of
logarithmic foliations.} Comp. Math., \textbf{142} (2006), 131--142.


\bibitem{GML}
X.  G\'{o}mez-Mont and A. Lins Neto, \emph{Structural stability of
singular holomorphic foliations having a meromorphic first
integral}. Topology  \textbf{30}  (1991),  no. 3, 315--334.


\bibitem{Griffiths-Harris}
P. Griffiths and J. Harris, \emph{Principles of Algebraic Geometry}, Wiley (1978).


\bibitem{Grauert}
H. Grauert and R.  Remmert, \emph{Theory of Stein spaces.}
Springer-Verlag (1979).


\bibitem{Jouanolou}
J. P. Jouanolou, \emph{\'{E}quations de Pfaff alg\'{e}briques.} Lecture
Notes in Mathematics, \textbf{708}. Springer, Berlin, 1979.


\bibitem{Malgrange}
B. Malgrange, \emph{Frobenius avec singularites - Le cas general},
Inventiones Math.  \textbf{39} (1977), 67-89.

\bibitem{Airton}
A. de Medeiros, \emph{Singular foliations and differential
$p$-forms.}  Ann. Fac. Sci. Toulouse Math. (6)  \textbf{9} (2000),
no. 3, 451--466.

\bibitem{coloquio}
J. V. Pereira,  \emph{Integrabilidade de Folhea\c{c}\~{o}es Holomorfas.}
Publica\c{c}\~{o}es Matem\'{a}ticas, IMPA,  2003.

\bibitem{richard}
R. W.  Richardson, \emph{A rigidity theorem for subalgebras of Lie
and associative algebras.} Illinois J. Math. \textbf{11} (1967)
92--110.


\bibitem{Saito}
K. Saito, \emph{On a generalization of de-Rham lemma.} Ann. Inst.
Fourier (Grenoble) \textbf{26} (1976), no. 2, vii, 165--170.


\end{thebibliography}
\end{document}